\documentclass[letterpaper, 10 pt, conference]{ieeeconf}   
\IEEEoverridecommandlockouts                              
\usepackage{graphics} 
\usepackage{epsfig} 
\usepackage{amsmath,amssymb,amsfonts}
\usepackage{graphicx}  
\usepackage{tikz} 
\usepackage{epstopdf}
\usepackage{cite}
\usepackage{algorithm,algorithmic,setspace}
\usepackage{url}
\usepackage{hyperref}
\hypersetup{
	colorlinks=true,
	linkcolor=blue,
	filecolor=blue,
	urlcolor=blue,
	citecolor = blue,
}

\newcommand{\mycomment}[1]{{\hfill\ttfamily\footnotesize\fontdimen2\font=0.3em\textcolor{blue}{/*#1*/}}}  
\renewcommand{\det}{\mathsf{det}}
\newcommand{\tr}{\mathrm{tr}}
\newcommand{\T}{^{\top}}  
\newcommand{\ie}{\textit{i.e.}} 
\newcommand{\SO}{\mathsf{SO}}
\newcommand{\SE}{\mathsf{SE}}
\newcommand{\0}{{0}}
\newcommand{\so}{\mathfrak{so}}
\newcommand{\se}{\mathfrak{se}} 
\newcommand{\Ad}{\mathsf{Ad}}

\allowdisplaybreaks

\title{\LARGE \bf 
	Observers Design for Inertial Navigation Systems: A Brief Tutorial 
}

\author{Miaomiao Wang and Abdelhamid Tayebi 
	\thanks{This work was supported by the National Sciences and Engineering Research Council of Canada (NSERC).}
	\thanks{The authors are with the Department of Electrical and Computer Engineering, Western University, London, ON N6A 3K7, Canada. A. Tayebi is also with the Department of Electrical Engineering, Lakehead University, Thunder Bay, ON P7B 5E1, Canada.
		{\tt\small mwang448@uwo.ca}, {\tt\small atayebi@lakeheadu.ca}}%
}

\def\BibTeX{{\rm B\kern-.05em{\sc i\kern-.025em b}\kern-.08em
		T\kern-.1667em\lower.7ex\hbox{E}\kern-.125emX}}
\begin{document}

	\maketitle
	\thispagestyle{empty}
	\pagestyle{empty}
	
	\begin{abstract}
		The design of navigation observers able to simultaneously estimate the position, linear velocity and orientation of a vehicle in a three-dimensional space is crucial in many robotics and aerospace applications. This problem was mainly dealt with using the extended Kalman filter and its variants which proved to be instrumental in many practical applications. Although practically efficient, the lack of strong stability guarantees of these algorithms motivated the emergence of a new class of geometric navigation observers relying on Riemannian geometry tools, leading to provable strong stability properties. The objective of this brief tutorial is to provide an overview of the existing estimation schemes, as well as some recently developed geometric nonlinear observers, for autonomous navigation systems relying on inertial measurement unit (IMU) and landmark measurements.
	\end{abstract}
	\section{Introduction}
	The algorithms used for the determination of the position, velocity and orientation of an object moving in a three-dimensional space, referred to as navigation observers, are instrumental in many applications, such as autonomous underwater, ground and aerial vehicles. Before diving into more details on navigation observers, it is worth mentioning a few historical facts related to the evolution of navigation observers over the years. The attitude estimation, which is a crucial component of navigation observers, is perhaps one of the most important problems in aerospace engineering that has generated many theoretical and technological advances over the years. In the early days of the flying era, the attitude was determined through the integration of the angular velocity provided by bulky mechanical Gyroscopic instruments. In the early 1960s, a different and relatively more reliable approach, relying on vector observations, has been introduced in the aerospace community. This approach consists in using measurements (in the vehicle frame) of some known vectors in the inertial frame. These vector observations can be obtained using different types of sensors such as accelerometers and magnetometers included in cheap and tiny micro-electromechanical systems (MEMS)-based IMU sensors used in low cost and small drone applications, or sophisticated sensors such as star trackers for satellites and expensive space missions. Several static solutions to the attitude determination problem have been proposed since the early 1960s (see, for instance, \cite{Wahba}, \cite{shuster1979}). Obviously, these methods, although simple, do not perform well in the presence of measurement noise, which motivated their reinforcement with Kalman-type capabilities leading to dynamic attitude estimation algorithms considered as the workhorse of the aerospace industry (See the survey paper  \cite{Crassidis2007survey}). Although successful in many aerospace applications, these techniques have to be taken with extra care as they rely on liberalizations (approximations) and heavy computations. More recently, geometric nonlinear observers, evolving on the Special Orthogonal group $\SO(3)$, have emerged and shown their ability in efficiently handling the attitude estimation problem \cite{mahony2008nonlinear}, providing almost global asymptotic stability guarantees way beyond the local stability of the Kalman filters. Due to the motion space topology, almost global asymptotic (exponential) stability is the strongest result one can achieve via time-invariant smooth observers. To strengthen the stability properties, several solutions have been recently proposed in the literature such that the hybrid synergistic approach initiated as a control problem in \cite{mayhew2011synergistic} and further improved and adapted for observers design in \cite{berkane2017hybrid}.\\
	From the implementation point of view, in low-cost applications relying on IMU measurements, most of the existing attitude estimation techniques rely on the assumption that the vehicles's linear acceleration is negligible which allows to rely on the accelerometer as a gravity sensor in the vehicle frame. In applications involving non-negligible linear accelerations, one can use the so-called velocity-aided attitude observers that rely on IMU and linear velocity measurements \cite{bonnabel2009invariant,hua2010attitude,roberts2011}. The linear velocity (with respect to the inertial frame) can be obtained, for instance, from a Global Positioning System (GPS). GPS-aided navigation algorithms have been extensively studied in the literature, especially by Fossen and Johansen and their collaborators, and many interesting papers have been published such as \cite{bryne2017nonlinear}.
	However, the implementation of these observers, in GPS-denied environments, is not straightforward. A typical solution, in this situation, consists in using vision-aided inertial navigation observers, relying on the measurements of the angular velocity, linear acceleration, and some landmark positions (known in the inertial frame). These observers allow for a simultaneous estimation of the attitude, position and linear velocity, as well as the gyro-bias (and sometimes even the accelerometer bias).\\
	The most popular approaches used to solve this problem, in the aerospace community, are the extended Kalman filter (EKF) and its variants such as the additive EKF (AEKF) and the multiplicative EKF (MEKF). Different versions of these algorithms exist depending on which type of attitude representation is used and how the estimation errors are generated. For instance, for the attitude estimation problem, Euler-angles based EKF has been proposed in \cite{farrell1970attitude} and unit-quaternion based AEKF and MEKF have been proposed in \cite{bar1985attitude,leffens1982kalman}. On the other hand, an interesting invariant EKF (IEKF) approach, with provable local stability guarantees, for the navigation problem has been proposed in \cite{barrau2017invariant} relying on the symmetry preserving observer techniques developed in \cite{bonnabel2008symmetry,bonnabel2009non}. In contrast with the standard EKF and its variants, the IEKF has a state trajectory independent error propagation; a very interesting feature resulting from how the estimation errors are generated. More recently, motivated by the lack of strong stability guarantees of the EKF-type observers, nonlinear observers have been proposed in the literature such as \cite{hua2018riccati,BerkaneECC2019,wang2020hybrid,wang2020nonlinear}.  These observers are endowed with provable strong stability properties and, although designed in a deterministic setting, they exhibit nice filtering properties with respect to measurement noise, and relatively low computational overhead compared to their EKF-type counterparts. 
	
	This paper intends to provide a brief overview of the existing state observers for the navigation problem (\textit{i.e.,} simultaneous estimation of the attitude, position, velocity) using IMU and landmark measurements.
	We should also point out that the references included in this paper are not exhaustive due to space limitation and also to the widespread of this topic among different fields (aerospace, marine applications, robotics, control systems \textit{etc.}).

	\section{Preliminary Material}\label{sec:preliminary}
	\subsection{Notations}
	The sets of real, non-negative real, and nonzero natural numbers are denoted by $\mathbb{R}$, $\mathbb{R}_{\geq 0}$  and $\mathbb{N}_{>0}$, respectively. We denote by $\mathbb{R}^n$ the $n$-dimensional Euclidean space. The Euclidean norm of a vector $x\in \mathbb{R}^n$ is defined as $\|x\| = \sqrt{x\T x}$. The $n$-by-$n$ identity matrix and zero matrix are denoted by $I_n$ and  $\0_n$. For a matrix $A\in \mathbb{R}^{n\times n}$, we define $\lambda_{M}^A$ and $\lambda_{m}^A$ as the maximum and minimum eigenvalues of $A$. By $\text{diag}(\cdot)$, we denote the diagonal matrix. We define $\hat{x}, \dot{x}, x^+$ and $x\T$ as the estimate, time-derivative, discrete update and transpose of the state $x$. We use the notation $\mathcal{N}(\mu,\Sigma)$ to denote the Gaussian noise with mean $\mu$ and covariance $\Sigma$.

	\subsection{Kinematics and measurements}
	The 3-dimensional Special Orthogonal group $\SO(3)$ is defined as 
	$\SO(3):=\{R\in \mathbb{R}^{3\times 3}| R R\T = R\T R=I_3, \det(R)=+1\}$, and its \textit{Lie algebra}  is given by 
	$\mathfrak{so}(3):=\{\Omega\in \mathbb{R}^{3\times 3}| \Omega = -\Omega\T\}$. 
	Let $\{\mathcal{I}\}$ be the inertial frame and $\{\mathcal{B}\}$ be a frame
	attached to the center of mass of a rigid body. Consider the rotation matrix  $R\in \SO(3)$ as the attitude of the body-fixed frame $\{\mathcal{B}\}$ with respect to the inertial frame $\{\mathcal{I}\}$. Define the vectors $p\in \mathbb{R}^3$ and $v\in \mathbb{R}^3$ as the position and linear velocity of the rigid body expressed in frame $\{\mathcal{I}\}$.   
	The kinematic model of a rigid body navigating in a 3-dimensional space is given by
	\begin{subequations}
		\begin{align}
			\dot{R} & = R \omega^\times \label{eqn:R}\\
			\dot{p} & = v \label{eqn:p}\\
			\dot{v} & =  g  + Ra\label{eqn:v}
		\end{align}
	\end{subequations}
	where $g\in \mathbb{R}^3$ denotes the gravity vector in frame $\{\mathcal{I}\}$, $\omega\in \mathbb{R}^3$ denotes the  angular velocity of the body-fixed frame expressed in frame $\{\mathcal{B}\}$, and  $a\in \mathbb{R}^3$ denotes the ``apparent acceleration" capturing all non-gravitational forces applied 
	to the rigid body expressed in frame $\{\mathcal{B}\}$. The map $ (\cdot) ^\times: \mathbb{R}^3 \to \so(3)$ projects any vector in $\mathbb{R}^3$ to $\so(3)$ (skew-symmetric matrix) such that $ x^\times y= x\times y$ for any $x,y \in \mathbb{R}^3$ with $\times$ denoting the vector cross-product on $\mathbb{R}^3$.
	
	We assume that the rigid body system is equipped with an inertial-vision system, which combines an IMU and an on-board vision system. The gyro and accelerometer measurements from IMU, denoted by $\omega_m$ and $a_m$, are modeled as follows:
	\begin{align}
		\omega_m &= \omega + b_\omega + n_\omega \label{eqn:omega_m}\\
		a_m &= a + b_a + n_a \label{eqn:a_m}
	\end{align}
	where $b_\omega, b_a$ denote the biases in the gyro and accelerometer measurements, $n_a,n_\omega$ are independent zero mean Gaussian noise components.  For the sake of presentation simplicity, and the purpose of this tutorial, the gyro and accelerometer biases are assumed to compensated using appropriate calibration techniques, see for instance \cite{tedaldi2014robust}. 
	
	We consider a family of $N$ landmark position measurements obtained from a vision system. Let $p_i \in \mathbb{R}^3$ be the position of the $i$-th landmark expressed in frame $\{\mathcal{I}\}$. The landmark measurements in the body-fixed frame $\{\mathcal{B}\}$ are modeled as
	\begin{equation}
		y_i = R\T(p_i-p) + n_{y_i}, \quad i=1,2,\cdots,N \label{eqn:output_y}
	\end{equation}
	where $n_{y_i}$ denotes a zero mean Gaussian noise. Note that the landmark position measurement $y_i$ can be obtained, for instance, from a stereo-vision system \cite{hartley2003multiple}.  
	
	\section{EKF-based Filters for Inertial Navigation}
	\subsection{Classical EKF}
	Consider the following nonlinear stochastic dynamic system: 
	\begin{align}
		\dot{x} &= f(x,u) +  G(x)  n_{x} \\
		y &= h(x) + N(x) n_{y}  
	\end{align}
	where $x\in \mathbb{R}^n$, $u\in \mathbb{R}^m$ and $y\in \mathbb{R}^p$ are the state, input and output, respectively;   $n_x \sim \mathcal{N}(0,V)$ and $n_y\sim \mathcal{N}(0,Q)$ are independent zero mean Gaussian noises;  $G(x)\in \mathbb{R}^{n \times n}$ and $N(x)\in \mathbb{R}^{p\times p}$ are matrix-valued functions of $x$. Since the inertial navigation system usually combines a high sampling-rate IMU sensor and a low sampling-rate vision system, a continuous-discrete version of the EKF is presented  in Algorithm \ref{alg:1}. The output measurements are assumed to be available at discrete time instants $\{t_k\}_{k\in \mathbb{N}_{>0}}$. The  subscript $t_k$ indicates the value of the variable at the time instant $t_k$. Necessary conditions for the local asymptotic convergence of the EKF are discussed in \cite{song1992extended,krener2003convergence}.
	
	\begin{algorithm}  [!ht]  \small
		\caption{Continuous-discrete EKF} \label{alg:1}
		\setstretch{1.0}
		\begin{algorithmic}[1]
			\renewcommand{\algorithmicrequire}{\textbf{Input:}}
			\renewcommand{\algorithmicensure}{\textbf{Output:}}
			\REQUIRE  $u(t)$ for all $t\geq 0$ and  $y(t_k)$ with $k\in \mathbb{N}_{> 0}$. 		
			\ENSURE  $\hat{x}(t)$ for all $t\geq 0$
			\FOR {$  k\geq 1 $}  	
			\WHILE{$t\in [t_{k-1},t_k]$}
			\STATE \(\dot{\hat{x}} =  f(\hat{x},u)\)  
			\STATE \(\dot{P}_t   = A_t P_t + P_t A_t\T  + G_tVG_t\T \) \\ 
			\mycomment{$A_t = \left.\frac{\partial f}{\partial x}\right|_{x=\hat{x}}$ and $G_t = G(\hat{x})$}
			\ENDWHILE
			\STATE \(K_k  = P_{t_k}C_{t_k}\T(C_{t_k}P_{t_k}C_{t_k}\T + N_{t_k} QN_{t_k}\T )^{-1}  \) \\
			\mycomment{$C_{t_k} = \left.\frac{\partial h}{\partial x}\right|_{x = \hat{x}_{t_k}}$ and $ N_{t_k} = N(\hat{x}_{t_k})$} 
			\STATE \(\hat{x}_{t_k}^+  =\hat{x}_{t_k} + K_k (y_{t_k} - h(\hat{x}_{t_k},u_{t_k})) \) 	
			\STATE \(P_{t_k}^+  = P_{t_k}-K_kC_{t_k}P_{t_k}\)  
			\ENDFOR
		\end{algorithmic}
	\end{algorithm} 
	
	Since the rotation $R\in \SO(3)$ is a 3-by-3 matrix,  the classical EKF cannot be directly applied to the state estimation for inertial navigation systems. A possible solution is to consider a lower-dimensional parameterization of the rotation $R$, for example the Euler angles representation \cite{farrell1970attitude}. However, all the three-parameter (minimal) representations of the attitude have singularity problems \cite{stuelpnagel1964parametrization}. An alternative globally non-singular (but non-unique) attitude representation is the unit-quaternion representation, which combines a scalar and a 3-dimensional vector (see \cite{shuster1993survey} for more details). A direct  application of the EKF for attitude estimation is the AEKF derived in \cite{bar1985attitude} considering the four components of the unit-quaternion vector as independent variables, and using the linear estimation error $q-\hat{q}$ which does not necessarily result in a unit-quaternion. 
	
	\subsection{MEKF for Inertial Navigation} \label{sec:MEKF}
	To avoid the issue of using an estimation error that does not preserve the motion space of the unit-quaternion, the MEKF proposed in \cite{leffens1982kalman} for attitude estimation  relies on a multiplicative quaternion estimation error $\tilde{q}=q\otimes \hat{q} ^{-1}$, where $\otimes$ is a special quaternion multiplication that preserves the quaternion group structure \cite{shuster1993survey}. Assuming that $\hat{q}$ is close enough to $q$, $\tilde{q}$ can be approximated by $ [\frac{1}{2}\tilde{\theta}\T, 1]\T$, where $\tilde{\theta} \in \mathbb{R}^3$ is generated by the multiplication of a small rotation angle with a unit vector corresponding to the direction of the rotation associated to $\tilde{q}$. There are many practical implementations of the EKF based on this multiplicative quaternion error, for instance \cite{mourikis2007multi,hesch2013consistency}.
	
	Similarly to the  multiplicative quaternion estimation error, one can use the multiplicative rotational estimation error as $\tilde{R}=R\hat{R}\T$ which can be approximated by  $I_3 +\tilde{\theta}^\times$ for small rotational discrepancies (\textit{i.e.,}  $\hat{R} \simeq R$). Taking the linear error for the position and velocity (\ie, $\tilde{p}=p-\hat{p}$ and $\tilde{v} = v-\hat{v}$), one obtains the estimation error vector $\tilde{x}=[\tilde{\theta}\T,\tilde{p}\T,\tilde{v}\T]\T\in \mathbb{R}^{9}$ and the output error $z = [z_1\T, \dots, z_N\T]\T$  with $z_i  = y_i  -  \hat{R}\T (\hat{p}-p_i)$. Then, the linearized model for the estimation error is given by
	\begin{align}
		\dot{\tilde{x}} &= \underbrace{\begin{bmatrix}
				\0_3 & \0_3 & \0_3\\
				\0_3 & \0_3 & I_3\\
				-(\hat{R} a_m)^\times  & \0_3 & \0_3 
		\end{bmatrix}}_{A_t} \tilde{x} - \underbrace{\begin{bmatrix}
				\hat{R}   & \0_3 \\
				\0_3 & \0_3\\
				\0_3 & \hat{R}  
		\end{bmatrix}}_{G_t} n_x  \label{eqn:defAG1} \\  
		z  &= \underbrace{\begin{bmatrix}
				\hat{R}\T(p_1-\hat{p})^\times & -\hat{R} \T & \0_3 \\ 
				\vdots & \vdots & \vdots \\
				\hat{R} \T (p_N-\hat{p})^\times & -\hat{R} \T & \0_3
		\end{bmatrix}}_{C_{t}} \tilde{x}   + n_y  \label{eqn:defC1}
	\end{align}
	where $n_x =[n_\omega\T, n_a\T]\T \sim \mathcal{N}(0,V)$ and $n_y = [n_{y_1}\T,\dots,n_{y_N}\T]\T  \sim \mathcal{N}(0,Q)$.
	A continuous-discrete version of the MEKF for inertial navigation is given in Algorithm \ref{alg:2}. The main advantage of the MEKF with respect to the AEKF is that the estimate $\hat{R}$ remains in $\SO(3)$ for all times.  From \eqref{eqn:defAG1} and \eqref{eqn:defC1}, the MEKF leads to matrices $A_t$ and $C_t$ that depend on the trajectory, which implies that the performance of the MEKF depends on the initial conditions, which may cause the estimation to diverge in some situations \cite{barrau2017invariant}.

	\begin{algorithm}[!ht]  \small
		\caption{MEKF for Inertial Navigation}\label{alg:2} 
		\setstretch{1.0}
		\begin{algorithmic}[1]
			\renewcommand{\algorithmicrequire}{\textbf{Input:}}
			\renewcommand{\algorithmicensure}{\textbf{Output:}}
			\REQUIRE  $\omega_m(t), a_m(t)$ for all $t\geq 0$, and  $y_{i}(t_k)$ with $k\in \mathbb{N}_{> 0}$ and $ i=1,2,\dots, N$. 
			\ENSURE  $\hat{R}(t),\hat{p}(t),\hat{v}(t)$ for all $t\geq 0$
			\FOR {$k \geq 1 $}  
			\WHILE{$t\in [t_{k-1},t_k]$}
			\STATE \(\dot{\hat{R}} =   \hat{R}   {\omega}_m^\times \)  
			\STATE \(\dot{\hat{p}} ~ =  \hat{v}   \)  
			\STATE \(\dot{\hat{v}} ~ =   g  + \hat{R}a_m  \)  
			\STATE \(\dot{P}_t   = A_t P_t + P_t A_t\T  + G_tVG_t\T \)  
			\mycomment{ $A_t$ and $G_t $ defined in \eqref{eqn:defAG1}}
			\ENDWHILE 
			\setstretch{1.3}
			\STATE 	\(K_k   = P_{t_k}C_{t_k}\T(C_{t_k}P_{t_k}C_{t_k}\T + Q)^{-1} \) with {$C_{t_k}$ defined in \eqref{eqn:defC1}} 
			\STATE $z = [z_1\T, \dots, z_N\T]\T$  with $z_i  = y_{i,t_k}  -  \hat{R}_{t_k}\T (p_i - \hat{p}_{t_k})$
			\STATE Obtain  $\delta\tilde{\theta},\delta\tilde{p}, \delta\tilde{v} \in \mathbb{R}^3$ from \( [\delta{\tilde{\theta}}\T,\delta{\tilde{p}}\T, \delta{\tilde{v}}\T]\T  = K_k z\)   
			\STATE \( \hat{R}_{t_k}^+ =   \exp((\delta{\tilde{\theta}})^\times) \hat{R}_{t_k} \) 
			\STATE \( \hat{p}_{t_k}^+ =\hat{p}_{t_k} + \delta{\tilde{p}}\) 
			\STATE \( \hat{v}_{t_k}^+ =\hat{v}_{t_k} + \delta{\tilde{v}}\) 
			\STATE \( P_{t_k}^+ = P_{t_k}-K_kC_{t_k}P_{t_k}\)  	
			\ENDFOR
		\end{algorithmic}
	\end{algorithm}  
	
	\subsection{IEKF for Inertial Navigation}\label{sec:IEKF}
	
	The  IEKF relies on geometric estimation errors (for all the state variables) instead of the linear position and velocity estimation errors used in the MEKF. 
	As shown in \cite{barrau2017invariant}, the inertial navigation state variables $R,p$ and $v$ can be grouped in a single element that belongs to the extended \textit{Special Euclidean group} of order 3, denoted by  $\SE_{2}(3)$, which is defined as
	$
	\SE_{2}(3) = \left\{X= \mathcal{T}(R,v,p) |  R\in \SO(3), p, v\in \mathbb{R}^3 \right\}
	$
	with the map $\mathcal{T}: \SO(3)\times \mathbb{R}^3 \times \mathbb{R}^3 \to \SE_2(3)$  defined by  
	$$ 
	\mathcal{T}(R,v,p)  =\begin{bmatrix}
		R & v & p\\
		0_{1\times 3} & 1 & 0 \\
		0_{1\times 3} & 0 & 1
	\end{bmatrix}  \in \mathbb{R}^{5\times 5}.
	$$
	Let $T_X\SE_2(3)$ be  the \textit{tangent space} of $\SE_2(3)$ at point $X$. The \textit{Lie algebra} of $\SE_2(3)$, denoted by $\se_2(3)$, is given by
	$$
	\se_2(3) = \left\{  \left.\begin{bmatrix}
		\omega^\times & \alpha & \nu\\
		0_{2\times 3} & 0_{2\times 1} & 0_{2\times 1}
	\end{bmatrix} \right| \omega,\alpha,\nu \in \mathbb{R}^3  \right\}  .
	$$
	More details about the Lie group $\SE_2(3)$ can be founded in \cite{barrau2017invariant,wang2020hybrid}. The kinematics (\ref{eqn:R})-(\ref{eqn:v}) with measurement noise can be rewritten in the following compact form:
	\begin{equation}
		\dot{X} = f(X,u) + X n_x^\vee  \label{eqn:dX}
	\end{equation}
	where $u=[\omega_m\T, a_m\T]\T \in \mathbb{R}^6, n_x= [n_\omega\T,n_a\T,0_{1\times 3}]\T\in \mathbb{R}^9$, and the nonlinear map $f: \SE_2(3) \times \mathbb{R}^6 \to T_X\SE_2(3)$ and the map $()^\vee: \mathbb{R}^9 \mapsto \se_2(3)$ are defined as
	\begin{align}  \setlength\arraycolsep{3pt} 
		f(X,u)=\begin{bmatrix}
			R\omega^\times& g+ Ra   & v \\
			0_{1\times 3} & 0 & 0\\
			0_{1\times 3} & 0 & 0
		\end{bmatrix}, 
		n_x^\vee = \begin{bmatrix}
			n_\omega^\times & n_a & 0_{3\times 1}\\
			0_{1\times 3} & 0  & 0 \\
			0_{1\times 3} & 0  & 0 
		\end{bmatrix}.
		\label{eqn:def_f}
	\end{align}
	Let $r_i= [p_i\T~0~1]\T \in \mathbb{R}^{5}$ for all $i=1,2,\cdots,N$ be the new inertial reference vectors with respect to the inertial frame $\{\mathcal{I}\}$, and from (\ref{eqn:output_y}) their measurements expressed in the body-fixed frame $\{\mathcal{B}\}$ are given by
	\begin{equation}
		b_i  = [y_i\T ~0~1]\T  =X^{-1} r_i, \quad i=1,2,\cdots,N. \label{eqn:output_X}
	\end{equation}
	For the IEKF,  the right-invariant error $\tilde{X} = \hat{X}X^{-1}= \mathcal{T}(\hat{R}R\T, \hat{v}-\hat{R}R\T v, \hat{p}-\hat{R}R\T p) $  and the output error $z = [z_1\T,\dots,z_N\T]$ with $z_i = [I_3, 0_{3\times 2}](\hat{X} b_i - r_i) = \hat{R}y_i + \hat{p} - p_i, \forall i=1,2,\dots,N$ are considered. Applying the first-order approximation $\tilde{X}  \approx I_5 + \tilde{x}^\vee$ and assuming $\|x\|\|n_x\|\ll \|n_x\|$ and $\|x\|\|n_y\|\ll \|n_y\|$, one obtains the following  linearized  models  \cite{barrau2017invariant}
	\begin{align}  
		\dot{\tilde{x}} &= \underbrace{\begin{bmatrix}
				\0_3 & \0_3 & \0_3\\
				g^\times & \0_3 & \0_3\\
				\0_3 & I_3 & \0_3
		\end{bmatrix}}_{A_t} \tilde{x}  - \underbrace{\begin{bmatrix}
				\hat{R} & \0_3 & \0_3\\
				\hat{v}^\times \hat{R} & \hat{R} & \0_3\\
				\hat{p}^\times \hat{R} & \0_3 & \hat{R}
		\end{bmatrix}}_{G_t} n_x   \label{eqn:defAG2}  \\
		z & = \underbrace{\begin{bmatrix}
				p_1^\times & \0_3 & -I_3\\ 
				\vdots & \vdots & \vdots \\
				p_N^\times & \0_3 & -I_3\\
		\end{bmatrix}}_{C_t} \tilde{x}- \underbrace{\begin{bmatrix}
				\hat{R} & \cdots & \0_3\\ 
				\vdots & \ddots & \vdots \\
				\0_3 & \cdots & \hat{R}\\
		\end{bmatrix}}_{N_t} n_y \label{eqn:defCN2}
	\end{align}
	with $n_x  \sim \mathcal{N}(0,V)$, $n_y = [n_{y_1}\T,\dots,n_{y_N}\T] \sim \mathcal{N}(0,Q)$.
	The right invariant estimation error $\tilde{X}  = \mathcal{T}(\hat{R}R\T, \hat{v}-\hat{R}R\T v, \hat{p}-\hat{R}R\T p)$ is different from the estimation error considered in the MEKF. This right invariant error leads a nice feature which consists in the fact that the matrices $A_t,C_t$ defined in \eqref{eqn:defAG2}  and \eqref{eqn:defCN2} are constant on a much bigger set of trajectories than equilibrium points as it is the case for the MEKF. The IEKF (more precisely the RIEKF) for inertial navigation proposed in \cite{barrau2017invariant} is summarized in Algorithm \ref{alg:3}.  As shown in \cite[Theorem 6]{barrau2017invariant}, if there exist three non-collinear measurable landmarks, the IEKF is locally asymptotically stable.

	\begin{algorithm}[!ht]   \small
		\caption{IEKF for Inertial Navigation}\label{alg:3}
		\setstretch{1.0}
		\begin{algorithmic}[1]
			\renewcommand{\algorithmicrequire}{\textbf{Input:}}
			\renewcommand{\algorithmicensure}{\textbf{Output:}}
			\REQUIRE $u(t)=[\omega_m(t)\T, a_m(t)\T, 0_{1\times 3}]\T  $ for all $t\geq 0$, and  $b_{i}(t_k)$ with $k\in \mathbb{N}_{> 0}$ and $ i=1,2,\dots, N$		
			\ENSURE  $\hat{X}(t)$ for all $t\geq 0$
			\FOR {$ k \geq 1 $}  	
			\WHILE{$t\in [t_{k-1},t_k]$}
			\STATE \(\dot{\hat{X}} =   f(\hat{X},u)\)
			\STATE \(\dot{P}_t = A_t P_t + P_tA_t\T + G_t V G_t\T\)  
			\mycomment{$A_t$ and $ G_t $ defined  in \eqref{eqn:defAG2}}	 	
			\ENDWHILE
			\STATE $z = [z_1\T,\dots,z_N\T]$ with $z_i= [I_3, 0_{3\times 2}](\hat{X}_{t_k}b_{i,t_k} - r_i)$  
			\STATE 	\(K_k   = P_{t_k}C_{t_k}\T(C_{t_k}P_{t_k}C_{t_k}\T + N_{t_k}QN_{t_k}\T)^{-1}\) 
			\mycomment{ $C_{t_k}$ and $N_{t_k}$ defined  in \eqref{eqn:defCN2} }
			\STATE \(\hat{X}_{t_k}^+  =  \exp(K_k z) \hat{X}_{t_k}\)
			\STATE \(P_{t_k}^+  = P_{t_k}-K_kC_{t_k}P_{t_k}\)		 
			\ENDFOR
		\end{algorithmic}
	\end{algorithm}  
	\section{Geometric Nonlinear Observers for Inertial Navigation} 
	
	In this section, we present some recently developed Riemannian geometry-based nonlinear observers for inertial navigation, relying on  IMU and  landmark measurements, with strong stability guarantees. We  first start with the nonlinear observer designed on the Lie group $\SE_{2}(3)$ using continuous IMU and landmark measurements. Then, these results are extended to the case where the landmark measurements are intermittently available at some discrete instants of time. 
	
	\subsection{Continuous Landmark Position  Measurements} \label{sec:cont_obs}
	\subsubsection{Nonlinear observer design}
	We assume that the IMU and landmark  measurements  are continuous and noise-free. Consider a set of scalar weights $k_i>0, i=1,2,\dots N$ such that $\sum_{i=1}^N k_i=1$, and  define $p_c  = \sum_{i=1}^{N} k_i p_i$ as the weighted center of landmarks in frame $\{\mathcal{I}\}$. Let us introduce the homogeneous transformation $\bar{r}_i = X_c^{-1}r_i$ with $X_c = \mathcal{T}(I_3, 0_{3\times 1}, p_c) $ and $r_i$ considered in \eqref{eqn:output_X}. 
	For later use, we define $ \bar{r}:= [\bar{r}_1 ~\bar{r}_2~\cdots ~\bar{r}_N]\in \mathbb{R}^{5\times N}$ and $b :=  [b_1 ~b_2~\cdots ~b_N]\in \mathbb{R}^{5\times N}$ with $b_i$ given in \eqref{eqn:output_X}.  
	
	Motivated by the geometric nonlinear observers on $\SE(3)$ proposed in \cite{vasconcelos2010nonlinear,hua2015gradient,khosravian2015observers,wang2019hybrid}, a continuous time-invariant observer on $\SE_2(3)$ is given as follows \cite{wang2020hybrid}:
	\begin{align} 	
		&\dot{\hat{X}} = f(\hat{X},u) - \Ad_{X_c}(\Delta) \hat{X}, \label{eqn:observer_smooth}\\
		&\Delta =
		- \mathbb{P}_\mathcal{K} ((\bar{r}-X_c^{-1}\hat{X}b) K_N \bar{r}\T)  
		\label{eqn:innovation_term}	
	\end{align}
	where $\hat{X}(0)\in \SE_2(3)$, $K_N = \text{diag}(k_1,\cdots,k_N)$, the map $f$ is defined in \eqref{eqn:def_f}, and the \textit{adjoint map} $\Ad: \SE_2(3) \times \mathfrak{se}_2(3) \to \mathfrak{se}_2(3)$ is given by $ \Ad_{X}U   := X U X^{-1}$ for any $X\in \SE_2(3), U\in \se_2(3)$. 
	The gain map $\mathbb{P}_\mathcal{K}: \mathbb{R}^{5\times 5} \to \se_2(3)$  with a set $\mathcal{K}  = (k_R, K_p, K_v)$ and $k_R\in \mathbb{R}, K_v,K_p\in \mathbb{R}^{3\times 3}$,  is defined as follows  
	\begin{equation} \setlength\arraycolsep{3.2pt}
		\mathbb{P}_{\mathcal{K}}\left(\begin{bmatrix}
			A_1 & a_2 & a_3 \\
			a_4\T & a_6 & a_7 \\
			a_5\T & a_8 & a_9
		\end{bmatrix}\right) = \begin{bmatrix}
			k_R\mathbb{P}_a(A_1) & K_v a_2 & K_p a_3\\
			0_{1\times 3}  & 0 & 0 \\
			0_{1\times 3}  & 0 & 0
		\end{bmatrix}  \label{eqn:project_PK}
	\end{equation}
	where   $A_1\in \mathbb{R}^{3\times 3}, a_2, \cdots, a_5 \in
	\mathbb{R}^3$, $a_6,\cdots,a_9\in \mathbb{R}$, and  $\mathbb{P}_a(\cdot)$  denoting the anti-symmetric projection of $A_1$ such that $\mathbb{P}_a(A_1) := (A_1-A_1\T)/2$.
	
	Note that observer \eqref{eqn:observer_smooth}  consists of two parts: the prediction term $f(\hat{X},u)$ relying on the IMU measurements, and the innovation term $\Delta$ designed in terms of the landmarks measurements. The introduction of the homogeneous transformation $X_c$ allows to achieve a nice decoupling property between the rotational and translational estimation error dynamics as pointed out in \cite{wang2017globally,wang2019hybrid,wang2020hybrid}. 
	
	On the other hand, observer \eqref{eqn:observer_smooth}  can be explicitly rewritten  on  $\SO(3) \times \mathbb{R}^3 \times \mathbb{R}^3$ as follows:
	\begin{align} 
		\begin{cases}
			\dot{\hat{R}} &= \hat{R}(\omega + k_R\hat{R}\T\sigma_R)^\times\\
			\dot{\hat{p}} &= \hat{v} +  k_R \sigma_R^\times (\hat{p} -p_c) +  K_p y     \\
			\dot{\hat{v}} &= g+ \hat{R}a +  k_R \sigma_R^\times\hat{v} + K_v y 
		\end{cases}
		\label{eqn:observer_Rpv}
	\end{align}  
	where $\hat{R}(0)\in \SO(3)$,  $\hat{p}(0),\hat{v}(0)\in \mathbb{R}^3$, and $\sigma_R$ and $y$ obtained from $\Ad_{X_c}(\Delta)$ are given as
	\begin{align} 
		\sigma_R &  = \frac{1}{2}\sum_{i=1}^N k_i (p_i - p_c)^\times (p_i - \hat{p}-\hat{R}y_i) \label{eqn:sigmaR} \\
		y    &=    \sum_{i=1}^N k_i(p_i - \hat{p}-\hat{R}y_i)  \label{eqn:sigmap}
	\end{align} 
	In view of \eqref{eqn:sigmaR} and \eqref{eqn:sigmap}, both $\sigma_R$ and $y$ have the  term $p_i - \hat{p}-\hat{R}y_i$ similar to the vector $z$ defined in Section \ref{sec:IEKF} for the IEKF. However, in observer \eqref{eqn:observer_Rpv} the innovation term $\sigma_R$ is used for the attitude estimation and $y$ is considered as the feedback for the position and velocity estimation.
	
	\subsubsection{Closed-loop system and gain design}
	Consider the state estimation error $X_c^{-1}X\hat{X}^{-1}X_c = \mathcal{T}(\tilde{R},\tilde{v},\tilde{p})$ with $\tilde{R}=R\hat{R}\T$, $\tilde{v} = v-\tilde{R}\hat{v}$ and $\tilde{p} = p-\tilde{R}\hat{p} - (I_3-\tilde{R})p_c$ as the estimation errors for the attitude, linear  velocity and position, respectively. If $p_c = 0_{3\times 1}$, the state estimation errors considered here are equivalent to the estimation errors considered in the IEKF.   From \eqref{eqn:output_y}, $y$ and $\sigma_R$ can be rewritten in terms of the estimation errors as $y=\tilde{R}\T \tilde{p}$ and $\sigma_R= \psi(M\tilde{R}) $  
	with $M = \sum_{i=1}^N k_i (p_i-p_c)(p_i-p_c)\T$ and the map $\psi(\cdot)$   defined as  $\psi(A)  =\frac{1}{2}[a_{32}-a_{23}, a_{13}-a_{31}, a_{21}-a_{12}]\T$ for any  $A=[a_{ij}]_{1\leq i,j\leq 3} \in \mathbb{R}^{3\times 3}$. 
	From \eqref{eqn:R}-\eqref{eqn:v} and \eqref{eqn:observer_Rpv}, one obtains the following closed-loop system 
	\begin{subequations} 
		\begin{align}
			\dot{\tilde{R}} &= \tilde{R} (-k_R \psi(M\tilde{R}))^\times   \label{eqn:closed-loop1-R}\\
			\dot{\tilde{p}} &=    \tilde{v}     - \tilde{R}K_p \tilde{R}\T \tilde{p} \label{eqn:closed-loop1-p}\\
			\dot{\tilde{v}}&=    (I_3-\tilde{R})g -    \tilde{R} K_v  \tilde{R}\T \tilde{p}  \label{eqn:closed-loop1-v}  .
		\end{align}
	\end{subequations}
	The introduction of  $p_c$ in   \eqref{eqn:sigmaR}-\eqref{eqn:sigmap} and the state estimation error,  leads to an interesting decoupling property for the closed-loop system, where the dynamics of $\tilde{R}$ are independent of $\tilde{p}$ and $\tilde{v}$ as shown in  \eqref{eqn:closed-loop1-R}. Suppose that there exist at least three non-collinear landmarks and a set of scalars $k_i>0, i=1,2,\dots,N$ such that $M$ is positive semi-definite with three distinct eigenvalues.  Then, one can show that the equilibrium point $\tilde{R}=I_3$ of \eqref{eqn:closed-loop1-R} is almost globally asymptotically stable \cite{mahony2008nonlinear}. Therefore, with the convergence of $\tilde{R}$, the convergence of $\tilde{p}$ and $\tilde{v}$ can be achieved using the following gain design approaches:
	
	\begin{itemize}
		\item \textit{Fixed-gain design:} 
		Let the constant $k_R> 0$ and  the matrices $K_p$ and $K_v$ be chosen as $K_p = k_p I_3$ and $K_v = k_v I_3$ with constant scalars $k_p,k_v>0$. Define the new state vector $\mathsf{x} = [\tilde{p}\T,\tilde{v}\T]\T \in \mathbb{R}^6$. From \eqref{eqn:closed-loop1-p}-\eqref{eqn:closed-loop1-v}, the dynamics of $\mathsf{x}$ are given by
		\begin{equation} 
			\dot{\mathsf{x}} = \underbrace{\begin{bmatrix}
					\0_3 & I_3\\
					\0_3 & \0_3
			\end{bmatrix}}_{A} \mathsf{x}  -K \underbrace{\begin{bmatrix}
					I_3 & \0_3
			\end{bmatrix}}_{C} \mathsf{x} + \delta_g \label{eqn:dx1}
		\end{equation}
		with $K = [k_p I_3, k_v I_3]\T$, $\delta_g = [0_{1\times 3}, g\T(I_3-\tilde{R})\T]\T $.  
		The dynamics of $\mathsf{x}$ can be seen as a linear time-invariant system with an additional disturbance term introduced by the gravity.  Obviously, the additional term $\delta_g \to 0_{6\times 1}$    as $\tilde{R} \to I_3$.
		Choosing $k_p,k_v>0$, one can easily verify that matrix $A-KC$ is Hurwitz.
		
		\item \textit{Variable-gain design:} 
		We consider $k_R>0$ as a constant gain and $K_p,~K_p$ as time-varying gain matrices.  Define the new state vector $\mathsf{x} = [\tilde{p}\T R,\tilde{v}\T R]\T$, whose dynamics from  \eqref{eqn:closed-loop1-p}-\eqref{eqn:closed-loop1-v} are given by
		\begin{equation} 
			\dot{\mathsf{x}} = \underbrace{\begin{bmatrix}
					-\omega^\times & I_3 \\
					\0_{3} & -\omega^\times
			\end{bmatrix}}_{A_t} \mathsf{x} - K \underbrace{\begin{bmatrix}
					I_3 & \0_{3}
			\end{bmatrix}}_{C_t} \mathsf{x} + \bar{\delta}_g \label{eqn:dx2}
		\end{equation}
		with $K= [\hat{R}\T K_p\T \hat{R}, \hat{R}\T K_v\T \hat{R}]\T$, $\bar{\delta}_g = [0_{1\times 3}, g\T(R-\hat{R})]\T $. One can also show that $\bar{\delta}_{g} \to 0_{6\times 1}$   as $\tilde{R} \to I_3$.  
		The gain $K$ is chosen as $K = PC_t\T Q_t$ with $P$ being the solution to the following continuous Riccati equation:  
		\begin{equation}
			\dot{P} = A_t P + PA_t\T  - PC_t\T  Q_t  C_t P + V_t \label{eqn:CRE} 
		\end{equation}
		where $P(0) \in \mathbb{R}^{6\times 6}$ is  a symmetric positive definite matrix. 
		Given the pair $(A_t,C_t)$ uniformly observable \cite[Lemma 3]{wang2020hybrid} and the  matrices $V_t\in \mathbb{R}^{6\times 6}$ and $Q_t\in \mathbb{R}^{3\times 3}$ uniformly positive definite and bounded, from \cite{bucy1967global,bucy1972riccati}  the solution $P$ to \eqref{eqn:CRE} is well defined on $\mathbb{R}_{\geq 0}$ and there exist positive constants $0< p_m\leq p_M < \infty$ such that $p_m I_6 \leq P  \leq p_M I_6$ for all $ t\geq 0$.
	\end{itemize} 
	Due to the topology of the Lie group $\SO(3)$, it is impossible to achieve robust and global stability results with smooth (or even discontinuous) state observers \cite{koditschek1989application,mayhew2011synergistic,mayhew2013synergistic}. Hence, the best stability result one can achieve with the continuous time-invariant observers, for inertial navigation, is almost global asymptotic stability. Motivated by the hybrid observers designed on the Lie groups $\SO(3)$ and $\SE(3)$ in \cite{wu2015globally,berkane2017hybrid,berkane2017attitude,wang2019hybrid},  observer \eqref{eqn:observer_smooth}-\eqref{eqn:innovation_term} can be hybridized by introducing a resetting mechanism that prevents the state variables from converging to or getting stuck at the undesired equilibrium points of the closed-loop system \cite{wang2020hybrid}. 
	%
	Applying the framework of hybrid dynamical systems \cite{goebel2009hybrid,goebel2012hybrid,teel2013lyapunov},  the equilibrium point $(\tilde{R}=I_3, \|\tilde{p}\|=\|\tilde{v}\|=0)$ for the hybrid version of the observer \eqref{eqn:observer_Rpv}-\eqref{eqn:sigmap} (with both fixed-gain and variable-gain approaches) is uniformly globally exponentially stable \cite[Theorem 1 and Theorem 2]{wang2020hybrid}.  Moreover, the hybrid observer with time-varying gains can be further extended to handle IMU biases with  global exponential stability  guarantees as shown in \cite[Theorem 5]{wang2020hybrid}.
	
	\subsection{Intermittent Landmark Position Measurements} \label{sec:disc_obs}
	Now, we assume that the IMU measurements are continuous and the landmark measurements are available at some discrete instants of time $\{t_k\}_{k\in \mathbb{N}_{>0}}$. This assumption is motivated by the fact that the IMU sampling rates are much higher than those of the vision system. Suppose that there exit constants $T_m, T_M$ such that $0<T_m \leq t_{k}-t_{k-1} \leq T_M < \infty$ with $t_0=0$ and $k=1,2,\dots,N$. Two types of nonlinear observers  for inertial navigation, extended from the continuous observer \eqref{eqn:observer_Rpv}, are presented using the framework of hybrid dynamical systems. 
	
	\subsubsection{Nonlinear observer with discrete attitude estimation} \label{sec:NOL-D}
	The estimated attitude, position and velocity are obtained via hybrid dynamics consisting of a continuous integration of the  kinematics using the IMU measurements  between two  consecutive landmark measurements  (\ie, $t\in[t_{k-1},t_{k}], k\in \mathbb{N}_{>0}$), and a discrete update upon the arrival of the landmark measurements (\ie, $t=t_k, k\in \mathbb{N}_{>0}$).  The hybrid nonlinear observer is given as follows: 
	\begin{align}
		\underbrace{
			\begin{array}{ll}
				\dot{\hat{R}} &= \hat{R}(\omega)^\times   \\ 
				\dot{\hat{p}} &=  \hat{v}   \\
				\dot{\hat{v}} &=   g  + \hat{R}a
		\end{array} }_{t \in [t_{k-1}, t_k],~ k\in \mathbb{N}_{>0} } 
		\underbrace{
			\begin{array}{ll}
				\hat{R}^+ &= R_\sigma \hat{R}  \\ 
				\hat{p}^+ &=   R_\sigma (\hat{p} - p_c +   \hat{R} K_p \hat{R}\T y) + p_c     \\
				\hat{v}^+ &= R_\sigma (\hat{v}  + \hat{R}  K_v \hat{R}\T  y)
		\end{array} }_{t \in \{t_k\},~ k\in \mathbb{N}_{>0}}
		\label{eqn:observer2}
	\end{align}
	with $\hat{R}(0)\in \SO(3)$,  $\hat{p}(0),\hat{v}(0)\in \mathbb{R}^3$, $K_v,K_p\in \mathbb{R}^{3\times 3}$, and $y$   defined in \eqref{eqn:sigmap}.  The   term $R_\sigma $ for the attitude estimation  is designed based on the Cayley’s map \cite{shuster1993survey} given as  
	$$
	R_\sigma = \frac{(1-\|\sigma\|^2)I_3 + 2\sigma \sigma\T + 2\sigma^\times }{1 + \|\sigma\|^2} \in \SO(3)
	$$
	where $\sigma = -2 k_R \sigma_R $ with $k_R>0$ and $\sigma_R$ defined in \eqref{eqn:sigmaR}.  
	Observer \eqref{eqn:observer2} shares a similar structure of prediction and update as the EKF. The hybrid dynamics for the attitude estimation in \eqref{eqn:observer2}  are adapted from \cite{berkane2019attitude} in the case of synchronous measurements and rewriting $\sigma$  as $\sigma  = k_R \sum_{i=1}^N k_i (\hat{R}R\T a_i)^\times a_i$ with $a_i = p_i-p_c$ for all $i=1,2,\dots,N$.

	Consider the same estimation errors as in Section \ref{sec:cont_obs}. From \eqref{eqn:R} and \eqref{eqn:observer2}, one obtains the following closed-loop system for the attitude estimation: 
	\begin{align}
		\begin{cases}
			\dot{\tilde{R}}  = \tilde{R}(0_{3\times 1} )^\times & t \in [t_{k-1}, t_k],~ k\in \mathbb{N}_{>0}\\
			\tilde{R}^+  =  \tilde{R} R_\sigma\T  & t \in \{t_k\},~ k\in \mathbb{N}_{>0}.
		\end{cases}
		\label{eqn:error2} 
	\end{align}	
	Suppose that there exist at least three non-collinear landmarks such that one can generate at least two non-collinear vectors among the vectors $a_i=p_i-p_c, i=1,2,\dots,N$. Choosing $k_R$ small enough such that $k_R ({\tr(M) - \lambda_{m}^M})<1$, the attitude estimation error $\tilde{R}$ converges exponentially to $I_3$ for any initial condition $|\tilde{R}(0)|_I< \epsilon$ with $0<\epsilon<1$ \cite[Theorem 1]{berkane2019attitude}. The design of the gain matrices $K_p,K_v$ and the convergence of the position and velocity estimation errors will be discussed later. 
	
	\subsubsection{Nonlinear observer with continuous attitude estimation}
	Instead of updating the attitude estimate at each arrival time of the landmark measurements, we will make use of an auxiliary variable $\eta\in \mathbb{R}^3$, which remains constant between two consecutive landmark measurements  and updates upon the arrival of the landmark measurements. The estimated attitude is obtained through a continuous integration of the attitude kinematics using the gyro measurements and the auxiliary variable $\eta$. The estimated position and velocity are obtained via hybrid dynamics consisting of a continuous integration of the translational dynamics using the accelerometer measurements and the auxiliary variable $\eta$ between two consecutive landmark measurements, and a discrete update upon the arrival of the landmark measurements. The hybrid nonlinear observer is given as follows:
	\begin{align}
		\underbrace{
			\begin{array}{ll}
				\dot{\hat{R}} &= \hat{R}(\omega + \hat{R}\T \eta )^\times   \\
				\dot{\eta}  &= 0_{3\times 1}  \\
				\dot{\hat{p}} &= \eta^\times (\hat{p}-p_c) + \hat{v}   \\
				\dot{\hat{v}} &= \eta^\times  \hat{v} + g  + \hat{R}a
		\end{array}}_{t \in [t_{k-1}, t_k],~ k\in \mathbb{N}_{>0} }
		\underbrace{
			\begin{array}{ll}
				\hat{R}^+ &= \hat{R}  \\
				\eta^+    & = k_R\sigma_R   \\
				\hat{p}^+ &=  \hat{p}  + \hat{R} K_p \hat{R}\T y     \\
				\hat{v}^+ &=  \hat{v}  + \hat{R}  K_v \hat{R}\T  y
		\end{array} }_{t \in \{t_k\},~ k\in \mathbb{N}_{>0}}
		\label{eqn:observer3}
	\end{align}
	with $\hat{R}(0)\in \SO(3), \hat{p}(0), \hat{v}(0),\eta(0)  \in \mathbb{R}^3$,   $k_R>0$, $K_v,K_p\in \mathbb{R}^{3\times 3}$, and  $\sigma_R$ and $y$ defined in (\ref{eqn:sigmaR}) and (\ref{eqn:sigmap}), respectively. 
	Contrary to the MEKF in Section \ref{sec:MEKF}, the IEKF in Section \ref{sec:IEKF} and the nonlinear observer \eqref{eqn:observer2} where the attitude is updated intermittently, the attitude estimate from observer \eqref{eqn:observer3} is continuous  thanks to the auxiliary variable $\eta$ which takes care of the jumps upon the arrival of the landmark measurements. This feature (\textit{i.e.,} continuous attitude estimates) is desirable in practice, when dealing with observer-controller implementations.
	
	Consider the same state estimation errors as in Section \ref{sec:cont_obs}. From \eqref{eqn:R} and \eqref{eqn:observer3}, one obtains the following closed-loop system for the attitude estimation:
	\begin{equation} 
		\underbrace{
			\begin{array}{ll}
				\dot{\tilde{R}} &= \tilde{R}(- \eta )^\times \\
				\dot{\eta}  &= 0_{3\times 1}   
		\end{array}}_{t\in[t_{k-1}, t_k],~ k\in \mathbb{N}_{>0}} 
		\underbrace{
			\begin{array}{ll}
				\tilde{R}^+ &= \tilde{R}   \\
				\eta^+    & = k_R\psi(M\tilde{R})  
		\end{array} }_{t \in \{t_k\},~ k\in \mathbb{N}_{>0} }
		\label{eqn:error3}
	\end{equation} 	
	Suppose that there exist at least three non-collinear landmarks such that matrix $\bar{M} = \frac{1}{2}(\tr(M)I_3 - M)$ is positive definite. Then, for any initial condition $|\tilde{R}(0)|_I<\epsilon  \sqrt{\varsigma_{M}}$ with $\sqrt{\varsigma_{M}}=\lambda_{m}^{\bar{M}}/\lambda_{M}^{\bar{M}}$ and $0<\epsilon<1$, there exist $\|\eta(0)\|$ and $k_R$ small enough such that $\tilde{R}$ converges exponentially to $I_3$ \cite[Theorem 3]{wang2020nonlinear}. The design of the matrices $K_p,K_v$ and the convergence of the position and velocity estimation errors will be discussed next. 
	
	\subsubsection{Closed-loop system and gain design} From \eqref{eqn:R}-\eqref{eqn:v},  both observers \eqref{eqn:observer2} and \eqref{eqn:observer3} share the same dynamics for the position and velocity estimation errors given as
	\begin{align}
		\underbrace{
			\begin{array}{ll} 
				\dot{\tilde{p}} &=   \tilde{v}   \\
				\dot{\tilde{v}} &=   (I-\tilde{R})g
		\end{array}}_{t\in[t_{k-1}, t_k], k\in \mathbb{N}_{>0}}
		\quad
		\underbrace{
			\begin{array}{ll} 
				\tilde{p}^+ &=  \tilde{p}  -   R K_p R\T \tilde{p}    \\
				\tilde{v}^+ &=   \tilde{v} -   R K_v R\T \tilde{p}
		\end{array}  }_{t \in \{t_k\}, k\in \mathbb{N}_{>0} }.
		\label{eqn:error4}
	\end{align}	 
	Similar to the Section \ref{sec:cont_obs}, both fixed-gain and variable-gain  design approaches are presented as follows:
	\begin{itemize}
		\item \textit{Fixed-gain design:} Choose  $K_p = k_p I_3$ and $K_v = k_v I_3$ with constant scalars $k_p,k_v>0$. Define the new state vector $\mathsf{x} = [\tilde{p}\T,\tilde{v}\T]\T$, whose dynamics derived from \eqref{eqn:error4} are given as
		\begin{equation}
			\begin{cases}
				\dot{\mathsf{x}}  ~=A \mathsf{x} + \delta_g & t\in[t_{k-1}, t_k],k\in \mathbb{N}_{>0} \\
				\mathsf{x}^+ = (I-KC) \mathsf{x} & t \in \{t_k\},k\in \mathbb{N}_{>0}
			\end{cases}  \label{eqn:hx1}
		\end{equation}
		where $\delta_g  = [0_{1\times 3},~ g\T (I_3-\tilde{R})\T ]\T $, $K  =[ k_p I_3, k_v I_3]\T$, and the matrices $A$ and $C$ are given by \eqref{eqn:dx1}.  
		The gain $K$ is chosen such that there exists a symmetric positive definite  matrix $P$ satisfying the following condition:
		\begin{equation}
			A_g\T\Phi(\tau)\T P \Phi(\tau) A_g - P <0, ~\forall \tau \in [T_m,T_M]  \label{eqn:defXi_P}
		\end{equation}		
		with $\Phi(\tau) = \exp(A \tau)$ and $A_g = (I-KC)$. One can show that for any $\mathsf{x}(0)\in \mathbb{R}^6$, the  state $\mathsf{x}$ converges exponentially to zero \cite[Theorem 3]{wang2020nonlinear}.

		\item \textit{Variable-gain design:} We consider $K_p$ and $K_v$ as time-varying   matrices. Define the new state vector $ {\mathsf{x}}  = [\tilde{p}\T R,\tilde{v}\T R]\T \in \mathbb{R}^6$, whose dynamics derived from \eqref{eqn:R} and \eqref{eqn:error4} are given as
		\begin{equation}
			\begin{cases}
				\dot{\mathsf{x}}   =A_t \mathsf{x} + \bar{\delta}_{g}  & t \in [t_{k-1},t_k],~k\in \mathbb{N}_{>0}\\
				\mathsf{x}^+  = (I-KC_t) \mathsf{x} & t \in \{t_k\},~k\in \mathbb{N}_{>0}
			\end{cases}  \label{eqn:hx3}
		\end{equation}
		where $\bar{\delta}_g = [0_{1\times 3}, g\T (R-\hat{R}) ]\T$, $K =[K_p\T,K_v\T]\T$, and the matrices $A_t$ and $C_t$ are given by \eqref{eqn:dx2}. 
		The gain $K$ is chosen as $ K = PC_t\T(C_tPC_t\T + Q_t)^{-1} $ with $P$ being the solution of the following continuous-discrete Riccati equations:
		\begin{subequations}
			\begin{align}
				\dot{P} ~~&= A_tP + PA_t\T + V_t, ~ t\in [t_{k-1},t_k] \label{eqn:CDRE-C}\\
				P^+ &= P-KC_t P,\qquad \quad   t \in \{t_k\} \label{eqn:CDRE-D}
			\end{align}
		\end{subequations}
		where $k\in \mathbb{N}_{>0}$, $P(0)$ is symmetric and positive definite. Given the pair $(A_t,C_t)$ uniformly observable \cite[Lemma 3]{wang2020hybrid} and   matrices $V_t\in \mathbb{R}^{6\times 6}$ and $Q_t\in \mathbb{R}^{3\times 3}$  uniformly positive definite and bounded, from \cite{deyst1968conditions,barrau2017invariant}  the solution $P$ to   \eqref{eqn:CDRE-C}-\eqref{eqn:CDRE-D} is well defined on $\mathbb{R}_{\geq 0}$ and there exist positive constants $0< p_m\leq p_M < \infty$ such that $p_m I_n \leq  P \leq p_M I_n$ for all $ t\geq 0$. Then, one can show that for any $\mathsf{x}(0)\in \mathbb{R}^6$, the   state $\mathsf{x}$  converges exponentially to zero \cite[Theorem 9]{wang2020nonlinear}.

	\end{itemize}
	The optimization problem \eqref{eqn:defXi_P}   can be solved using either the polytopic embedding technique proposed in \cite{ferrante2016state} or the finite-dimensional LMI approach proposed in   \cite{sferlazza2019time}. 
	A complete procedure for solving the infinite-dimensional optimization problem \eqref{eqn:defXi_P}, adapted from \cite{sferlazza2019time}, is provided in \cite[Appendix A]{wang2020nonlinear}. 
	On the other hand, in order to tune the matrices $V_t$ and $Q_t$ using the covariance of the measurement noise, we consider the measurement noise as in \eqref{eqn:omega_m}-\eqref{eqn:output_y}.  Assuming $\|\mathsf{x}\|\|n_\omega\|\approx 0$,  the continuous dynamics of $\mathsf{x}$ in  \eqref{eqn:hx3} and  $y$ in \eqref{eqn:sigmap} can be approximately rewritten around $\|\mathsf{x}\|=0$ as
	\begin{align} 
		\begin{cases}
			\dot{\mathsf{x}} ~ = \underbrace{\begin{bmatrix}
					-\omega_m^\times & \0_3\\
					\0_3 & -\omega_m^\times 
			\end{bmatrix}}_{\bar{A}_t} \mathsf{x}  - \underbrace{\begin{bmatrix}
					(\hat{R}\T \bar{\hat{p}}) ^\times & \0_{3} \\
					(\hat{R}\T \hat{v})^\times & I_3
			\end{bmatrix}}_{G_t} u + \bar{\delta}_g \\
			y  ~= \bar{C}_t \mathsf{x} - \hat{R}\sum_{i=1}^N k_i n_y^i
		\end{cases}
		\label{eqn:hx3-new}
	\end{align}
	with $\bar{\hat{p}} = \hat{p}-p_c$, $u=[n_\omega\T,n_a\T]\T \sim \mathcal{N}(0,\bar{V})$ and $n_{y_i}\sim   \mathcal{N}(0,\bar{Q}), \forall i=1,2,\dots,N$. 
	Then, the matrices $V_t$ and $Q_t$ can be chosen as  
	\begin{align} \textstyle
		V_t~  = G_t  \bar{V}  G_t\T +  \varepsilon I_6, \quad  
		Q_t  = \sum_{i=1}^N k_i^2 \hat{R} \bar{Q}\hat{R}\T \label{eqn:VQ1}
	\end{align}
	where a small scalar $\varepsilon>0$ is introduced to ensure that $V_t$ is uniformly positive definite, and we made the assumption that the noise signals in the landmark measurements are uncorrelated. 
	The implementation procedure of the nonlinear observer \eqref{eqn:observer3} with variable gains is summarized in Algorithm \ref{alg:4}. 
	As shown in \cite{wang2020nonlinear}, the computational costs of the nonlinear observers presented in this subsection are lower than those of the IEKF. 
	Note that to achieve global exponential stability as shown in \cite{wang2020hybrid} one can hybridize the observer by introducing some appropriately designed jump mechanisms on the estimates $\hat{R},\hat{p}$ and $\hat{v}$ as in \cite{wang2020hybrid}.

	\begin{algorithm}[!ht]  \small
		\caption{Nonlinear   Observer for Inertial Navigation} \label{alg:4}
		\setstretch{1.0}
		\begin{algorithmic}[1] 
			\renewcommand{\algorithmicrequire}{\textbf{Input:}}
			\renewcommand{\algorithmicensure}{\textbf{Output:}}
			\REQUIRE  $\omega_m(t), a_m(t)$ for all $t\geq 0$, and  $y_{i}(t_k)$ with $k\in \mathbb{N}_{> 0}$ and $ i=1,2,\dots, N$. 
			\ENSURE  $\hat{R}(t),\hat{p}(t),\hat{v}(t)$ for all $t\geq 0$
			\FOR {$  k \geq 1$}   
			\WHILE{$t\in [t_{k-1},t_k]$}
			\STATE \(\dot{\hat{R}}  = \hat{R}(\omega_m  + \hat{R}\T \eta)^\times  \)  
			\STATE \(\dot{\eta}  = 0_{3\times 1}  \)  
			\STATE \(\dot{\hat{p}}  = \eta^\times (\hat{p}-p_c) + \hat{v} \)  
			\STATE \(\dot{\hat{v}}  =  \eta^\times \hat{v} + g  + \hat{R}a_m  \)  
			\STATE \(\dot{P}_t  = \bar{A}_t P_t + P_t\bar{A}_t\T + V_t \)  
			\mycomment{$\bar{A}_t$ defined  in \eqref{eqn:hx3-new} and $V_t$ defined  in \eqref{eqn:VQ1}}
			\ENDWHILE  
			\STATE $\sigma_R    = \frac{1}{2}\sum_{i=1}^N k_i (p_i - p_c)^\times (p_i - \hat{p}_{t_k}-\hat{R}_{t_k}y_{i,t_k}) $  
			\STATE $	y     =    \sum_{i=1}^N k_i(p_i - \hat{p}_{t_k}-\hat{R}_{t_k}y_{i,t_k}) $   
			\STATE 	\(K_k   = P_{t_k} {C}_{t_k}\T( {C}_{t_k}P_{t_k} {C}_{t_k}\T + Q_{t_k})^{-1} \) \mycomment{$C_{t_k}$ defined in \eqref{eqn:dx2}  and $Q_{t_k}$ defined  in \eqref{eqn:VQ1} }
			\STATE Obtain $K_p,K_v$ from \( [K_p\T  ,  K_v\T  ]\T = K_k\)  
			\setstretch{1.3}
			\STATE \( \hat{R}_{t_k}^+ =\hat{R}_{t_k}\) 
			\STATE \(\eta_{t_k}^+     = k_R \sigma_R\) 
			\STATE \(\hat{p}_{t_k}^+ = \hat{p}_{t_k} +   \hat{R}_{t_k}K_p \hat{R}_{t_k}\T  y  \) 
			\STATE \( \hat{v}_{t_k}^+ = \hat{v}_{t_k} +    \hat{R}_{t_k}K_v \hat{R}_{t_k}\T  y \) 
			\STATE \( P_{t_k}^+ = P_{t_k}-K_kC_{t_k}P_{t_k}\)  	
			\ENDFOR		
		\end{algorithmic}
	\end{algorithm} 
	\section{Conclusion}
	In this brief tutorial, we provide an overview of the most popular estimation techniques, as well as some recently developed ones, for the autonomous navigation problem using IMU and landmark measurements. We tried to provide the reader with practical implementation algorithms leaving aside the technical details that may be found in the corresponding papers. Several challenging problems remain open in this field, such as the design of nonlinear navigation observers, endowed with global asymptotic (exponential) stability, using biased IMU and monocular vision systems (with direct pixel measurements or bearing measurements). 
	Another interesting and challenging problem consists in the design of monocular-vision based nonlinear geometric Simultaneous Localization and Mapping (SLAM) observers, endowed with strong stability guarantees.
	
	\bibliographystyle{IEEEtran}
	\bibliography{mybib}
	

\end{document}